\documentclass[10pt]{article}

\usepackage{amsmath,amssymb}

\newcommand{\comment}[1]{}

\newcommand{\cov}{\vartriangleleft}

\newcommand{\primconst}[1]
        {\mbox{\sf #1}}
\newcommand{\Pos}{\primconst{Pos}}

\newtheorem{thm}[subsection]{Theorem}
\newtheorem{lem}[subsection]{Lemma}
\newtheorem{prop}[subsection]{Proposition}
\newtheorem{cor}[subsection]{Corollary}

\newtheorem{rmk}[subsection]{Remark}
\newenvironment{remark}%
      {\begin{rmk}\rm}
      {\end{rmk}}

\newcommand{\proof}{{\noindent \bf Proof. }}

\newcommand{\enproof}{$\Box$\\}

\newcommand{\Pow}{\textsf{Pow}}

\begin{document}
\title{On some peculiar aspects of the constructive theory of point-free
spaces\footnote{To appear on: Math. Log. Quart. 56, No. 4, 1 – (2010) / DOI 10.1002/malq.200910014. }}

\author{Giovanni Curi \\ \\
{\small Dipartimento di Informatica - Universit\`a di Verona}\\
{ \small  Strada le Grazie 15 - 37134 Verona, Italy.}\\
{ \small e-mail: giovanni.curi@univr.it}}

\date{}
\maketitle

\begin{abstract}
This paper presents several independence results concerning the  topos-valid and the intuitionistic (generalised) predicative theory of locales. In particular, certain consequences of the consistency of a general form of Troelstra's uniformity principle with constructive set theory and type theory are examined.
\end{abstract}

\noindent \textit{Key words}: Locales, formal spaces, constructive set theory and type theory, topos logic,  independence results, uniformity principle.

\medskip

\noindent \textit{MSC (2010)}: 03F50, 06D22, 54A35.

\section{{\normalsize Introduction}}
It may be argued that the well-known equivalence of theorems such as Tychonoff theorem, or  Stone-\v{C}ech compactification, with the axiom of choice, or other similarly non-constructive principles, is a consequence of the chosen formulation of the concept of space, rather than being an intrinsic feature of these results.
By replacing the ordinary notion of topological space with that of \textit{locale} (or frame, or complete Heyting algebra) one obtains fully general versions of these theorems  that can be proved without any choice, and often with no application of the principle of excluded middle  \cite{J83,FS79,FG82,J82}.

The notion of locale is for this reason the concept of space generally adopted in  choice-free intuitionistic settings, such as toposes or intuitionistic set theory  (IZF) \cite{MM}.
By not assuming  as available impredicative principles as the existence of powersets,
the concept of \emph{formal space}, or \emph{set-generated locale}, plays a corresponding role in even weaker systems,   as constructive set theory (CZF) or constructive type theory (CTT)
\cite{Sa87,CSSV,A,Pa,CuSC}. The main criterion of adequacy  of this notion is that,
considered in fully impredicative settings as ZF or IZF\footnote{In this paper a system is defined    `fully impredicative' if it is at least as strong as full higher order arithmetic HHA (topos logic).},
the category FSp  of formal spaces is equivalent to the ordinary category of locales.

This paper deals with certain peculiar features of both the constructive and the intuitionistic theory of locales. In fact, we will mainly be concerned with certain  independence results that follow from  the consistency of CZF and CTT with a generalised form of Troelstra's \textit{principle of uniformity} \cite{TvD}.

Our main results are mostly related to the following fundamental `structural' aspect of the theory of locales: considered in any topos,  the category of locales is complete and cocomplete, i.e., all limits (in particular products) and colimits exist in this category.
By contrast,   the existence of  binary products of arbitrary formal spaces already seems to require the use of strongly impredicative principles, that are not available in the generalised predicative settings under consideration.
In particular to remedy this deficiency,  the concept of \emph{inductively generated formal space} was introduced in \cite{CSSV,A}: inductively generated formal spaces define a full subcategory FSp$_i$ of the category  of formal spaces, in which limits and colimits do exist (albeit under the assumption of strong principles for the existence of inductively defined sets, such as the axiom REA in constructive set theory, see e.g. \cite{Pa,A}).

By exhibiting a particular formal space that CTT cannot prove to be inductively generated,  FSp$_i$  has been shown  to form a proper subcategory of FSp   \cite{CSSV} (see \cite{G} for a similar result in CZF).
Nevertheless, since FSp$_i$ contains in particular all locally compact formal spaces, and
since, considered in a fully impredicative setting, this category is still equivalent to the category of locales (as every formal space is inductively generated  in such a setting), the concept  of inductively generated formal space has generally been regarded as providing the proper constructive analogue of the notion of locale.

In this paper we show  that the restriction to the category FSp$_i$ is, however, a very severe one:
we prove that CTT, CZF, as  several extensions of CZF, including REA and the impredicative unbounded separation scheme, cannot prove that a non-trivial \textit{Boolean formal space} - i.e., a formal space whose associated frame is a non-trivial complete Boolean algebra -, is inductively generated.
This result provides us with an example, for every given formal space $S$ (inductively generated or not), of a formal space that these systems cannot prove to be inductively generated, namely the least dense subspace of $S$.
Similar facts also hold for De Morgan (or extremally disconnected) formal spaces, and for formal spaces whose associated frame is the the Dedekind--MacNeille completion of a poset.

Further independence results, concerning compactness, overtness (openness) and existence of points, will then be shown to hold with respect to the internal language of toposes (HHA), IZF, and/or CTT and CZF.
In particular, we show that
CZF (+REA+...), CTT cannot prove that a non-trivial formal space is compact and De Morgan. This is in contrast with a well-known result of M. H. Stone,  valid in any topos: in HHA, or IZF,
the frame of ideals of a complete Boolean algebra is a compact De Morgan locale \cite{Johnstone80}.
It follows in particular that, for no non-trivial compact regular formal space $S$, the  Gleason cover of  $S$ \cite{Johnstone80,J82} can be constructed in CZF, CTT.

The paper is organized as follows: basic facts on formal spaces/locales as treated in constructive settings are recalled in Section \ref{Preliminaries}.  In Section \ref{NFS} the  version of the uniformity principle that we shall exploit is presented and its incompatibility with De Morgan law is exhibited.
The independence results  concerning Boolean and De Morgan formal spaces are described in Section \ref{BFS}; the case of spaces arising via the Dedekind--MacNeille completion of a poset,  and a problem left open  in \cite{CSSV}, are discussed in Section \ref{OtherEx}.

\section{{\normalsize Preliminaries}}\label{Preliminaries}
The reader is referred to  \cite{AR,NPS90} for background on Aczel's
constructive set theory (CZF) and constructive type theory (CTT), respectively.
In the following,
we shall use CZF strictly to indicate the basic formulation of Aczel's theory.
An  extension of CZF that is often considered, particularly in connection with constructive
locale theory, is the theory CZF+REA. The regular extension axiom REA is needed to ensure that certain inductively defined classes  are sets \cite{AR}. Extending CZF with the full separation scheme (Sep) and the powerset axiom yields the fully impredicative set theory IZF. Adding the law of excluded middle to either theory (CZF or IZF) gives ZF.
In the following, we use HHA (for intuitionistic higher-order Heyting arithmetic)
to indicate topos logic \cite{MM,TvD}.

General information on locales
may be found in \cite{FG82,FS79,J82}; for basic
facts concerning the theory of formal spaces in constructive
predicative settings such as CZF, CTT see
\cite{A, CSSV, Sa87,Pa} or \cite{AczelCuri,CuSC}. Here we
synthetically recall the notions needed in this note.

A \textit{formal topology}, or \emph{formal space}, is a pair ${S}
\equiv (S,\cov)$ where $S$, the \textit{base}, is a set, and $\cov$,
the \emph{covering relation}, is a relation between elements and
subsets of $S$ satisfying:

\begin{itemize}
  \item  [$i.$] $a\in U$ implies $a\cov U$,
  \item  [$ii.$] if $ a\cov U$ and $U\cov V$, then $a\cov V$,
  \item [$iii.$] $a\cov U\mbox{ and } a\cov V$ imply $a\cov U\downarrow
  V$,
\end{itemize}

\noindent where
$ U \cov V \equiv (\forall u \in U)u \cov V,$ and $U\downarrow
V\equiv \{ d \in S :\ (\exists u \in U)\ (d \cov \{ u \})\ \&\
(\exists v \in V)\ (d \cov \{ v \}) \}.$
In CZF, the covering is formally a subclass of $S\times \Pow(S)$, where $\Pow(S)$ is the \emph{class} of subsets of $S$; in addition to  $i-iii$,  a further requirement in that context is that the class
${\cal S}(U)\equiv\{a:a\cov U\}$ be a set for all $U$ (see \cite{A}
for more). Two subsets $U,V$ of $S$ are the same \textit{formal
open}, $U=_S V$, exactly when  $U\cov V \ \& \ V\cov U$. Observe
that one may always assume that $S$ has a `top' element, i.e., an
element $1_S$ such that $S=_S \{1_S\}$. An implication operation is defined on formal opens by   $U\to V\equiv \{a\in S: a\downarrow U\cov V\}$. The \emph{pseudocomplement}
$U^*$ of $U$ (the largest open disjoint from $U$) is given by:
\[U^*\equiv U\to \emptyset\equiv \{a\in S:a\downarrow U\cov \emptyset\}.\]

A \emph{morphism} $f:S_1\longrightarrow S_2$ of formal topologies is  a mapping
 $f:S_1\longrightarrow \Pow(S_2)$  satisfying, for all $a,b\in S_1$, $U\subseteq S_1$,

\begin{itemize}
\item[$i.$] $f(S_1) =_{S_2} S_2$,
\item[$ii.$] $f(a)\downarrow  f(b)\cov  f(a\downarrow  b)$,
\item[$iii.$] $a \cov U$ implies $f(a) \cov f(U)$
\end{itemize}

\noindent where, for $U$ a subset of $S_1$, $f(U)\equiv\bigcup_{a\in U}
f(a)$. Two morphisms $f,g:S_1 \longrightarrow S_2$ are defined to be equal precisely when
$f(a)=_{S_2} g(a)$ for all $a\in S_1$. On any formal topology $S$,
the identity morphism is given by $id_S(a)=\{a\}$, for all $a$.

\medskip

A \textit{(formal) point} of a formal space $S$ is  a subset
$\alpha$ of $S$  satisfying:

\begin{itemize}
  \item [$i.$] $(\exists a\in S) a\in \alpha$
  \item [$ii.$] $a,b\in \alpha$ implies $(\exists c\in \alpha)$  $c\in
a\downarrow b$.
  \item [$iii.$] $a\in \alpha$ and $a\cov U$ imply $(\exists b\in U)$
$b\in \alpha$.
\end{itemize}

\noindent In particular, the top element $1_S$ (when it is present) belongs
to every point, and for no $a\in \alpha$, is $a\cov \emptyset$.

Even classically, a formal space may well have no points and be
\textit{non-trivial}, i.e., such that $\neg (S\cov \emptyset)$. In terms of logic, this is because points of a formal space
are the $\Pow(\{\top\})$-valued (classically two-valued) models of a
geometric theory, which may be consistent without having a model
\cite{FG82}.

\medskip

A \textit{subspace}  of a formal space ${S} \equiv(S,\cov)$ is a
formal space ${S}'\equiv(S,\cov')$, on the same base, and with
$\cov'$ satisfying $i. \cov\subseteq \cov'$, and $ii.\ x\downarrow '
y \cov'  x\downarrow y$. See  \cite{CuSC} for a more detailed discussion. For example, a
(formal) open subset $V\subseteq S$ determines the closed subspace
$S^V\equiv(S,\cov^V)$, with $a\cov^V U \iff a\cov U\cup V$
(intuitively, $S^V$ represents the complement of the open $V$ as a
subspace).

\medskip

A formal space $S$ is \textit{set-presented} iff there are families of sets $I(x)$, for $x$ in $S$, and $C(x,i)\subseteq S$, for $x\in S$,  $i\in I(x)$, such that \[a\cov U \iff (\exists i\in I(a)) C(a,i)\subseteq U.\]
Observe that
this implies $a\cov C(a,i)$ for all $i$. In CZF+REA, CTT,  $S$ is set-presented if and only if it is \emph{inductively generated} in the sense of \cite{A,CSSV}.
In CZF, CTT, the class of
set-presented formal spaces contains all locally compact spaces
\cite{A}. In a topos, or in (I)ZF, all formal spaces are
trivially set-presented: one simply defines $I(x)=\{U\in
\Pow(S):x\cov U\}$, $C(x,U)=U$.
The full subcategory of set-presented formal spaces has limits and colimits in sufficiently strong versions of constructive set theory and type theory.

In CZF, a \emph{class-frame} (or class-locale) $L$ is a partially
ordered class that has a top element, binary meets, and suprema for
arbitrary \emph{sets} of elements of $L$, and that is such that meets distribute over the set suprema.
A class-frame is said to be
\emph{set-generated} by a subclass $B$ if: $i.$ $B$ is a set; $ii.$ the class $\{b\in B:b\leq x\}$ is a set and $x=\bigvee\{b\in B:b\leq
x\}$, for all  $x\in L$.

\noindent Morphisms of set-generated frames are class-functions
respecting meets, the top, and arbitrary set joins.

Given a formal topology $S$, let  the
collection of \emph{saturated subsets} of $S$, i.e., the class
$\{U\subseteq S: {\cal S}(U)=U\}$, be denoted by $Sat(S)$. Endowed with the operations $U
\wedge V\equiv U \cap V=U\downarrow V\ \textrm{and}  \ \bigvee
_{i\in I} U_{i}\equiv {\cal S}(\bigcup_{i\in I}U_i)$, $Sat(S)$ is a
set-generated frame. The implication  operation previously recalled defines an implication  operation on $Sat(S)$ in the usual sense, making it in a complete Heyting algebra. In particular, $U^*$, for $U\in Sat(S)$ is the pseudocomplement of $U$ in the ordinary lattice-theoretic sense.

With their respective morphisms, formal topologies and set-generated
class-frames form equivalent categories  \cite{A, AczelCuri}. With
powersets, every set-generated class-frame has a set of elements, so
it is just an ordinary frame (locale). Therefore, in fully
impredicative settings such as toposes, the category FT  of
formal topologies is equivalent to that of frames (see also
\cite{Sa87}). Its opposite FSp=FT$^{op}$, here referred to
as the category of \textit{formal spaces} (often simply spaces) and
{\it continuous functions}, is thus  equivalent in such settings to
the category of locales.

\section{{\normalsize Uniformity principles}}\label{NFS}
To distinguish the behavior of formal spaces in constructive
settings from that in an intuitionistic but fully impredicative
context, we will exploit a generalised form of the so called uniformity principle \cite{TvD,Ro}.
In constructive set theory this is so formulated: for every set $I$,
\[(\forall x)(\exists y\in I)A(x,y)\to
(\exists y\in I)(\forall x)A(x,y)\ \ \ \textbf{(GUP-CZF)}.\]

In \cite{BennoThesis,BergMoerdijk} this principle has been proved to be consistent (in particular)
with CZF+REA+ PA+Sep, where REA is the regular extension axiom, PA
is the presentation axiom, and Sep is impredicative unbounded
separation (see \cite{BergMoerdijk} for a list of other principles compatible with GUP-CZF. Consistency of these principles with CZF is shown in \cite{BergMoerdijk} by the definition of a model that has independently been noted also in \cite{Lubarsky2006} and  \cite{Streicher}; see also \cite{Rathjen06}).
Note that GUP-CZF follows from its instance:
\[(\forall x)(\exists y\in \omega)A(x,y)\to
(\exists y\in \omega)(\forall x)A(x,y)\ \ \ \textbf{(UP-CZF)}\]
($\omega$ is the set of natural numbers), and the principle that every set is subcountable,  also  valid in the model of GUP described in \cite{BergMoerdijk, Lubarsky2006, Streicher}.
It will be convenient to note explicitly the following consequence of GUP-CZF: for every set $I$,
\[ (\forall p\in
\Pow(\{\top\}))(\exists i\in I)A(p,i)\to (\exists i\in I)(\forall
p\in \Pow(\{\top\}))A(p,i)\ \ \ \textbf{(GUP$'$-CZF)}\]
\noindent where
$\Pow(\{\top\})$ is the powerclass of the one-element set
 (the antecedent of GUP$'$-CZF yields
$(\forall x)(\exists y\in I)((\exists z)((\forall w) (w\in z \leftrightarrow w\in x \ \& \ w\in \{\top\})\ \& \ A(z,y)))$; one can then apply GUP-CZF).
The type-theoretic formulation of this principle, first exploited in \cite{CSSV}, is recalled in the Appendix.

\medskip

We write \textbf{EM}, \textbf{DML} for the principle of excluded middle and De Morgan law, respectively.
Recall that De Morgan law $ \neg (P \wedge Q)\to \neg P \vee \neg Q$ for all propositions $P,Q$,
is equivalent to \[\neg \neg P \vee \neg P\] for all $P$.
By the identification of subsets of the one-element set with
restricted formulas  (those in which all quantifiers are
bounded) \cite{AR}, in CZF  this principle for restricted
formulas can be formulated as \[(\forall p \in \Pow(\{\top\}))
p^{**}\cup p^*=\{\top\}\ \ \ (\textbf{[R]DML})\] where $p^*\equiv \{x\in \{\top\}:x\not\in p\}$. Note that,
considered in IZF, [R]DML  expresses De Morgan Law for arbitrary formulas.

The generalised uniformity principle conflicts in CZF with [R]DML, and with DML in CTT.
We prove the first fact: assume that $p^{**}\cup p^*=\{\top\}$ for all
$p\in \Pow(\{\top\})$. Define a relation $F\subseteq
\Pow(\{\top\})\times \{0,1\}$ by letting \[(x,y)\in F\iff (\top\in
x^{**}\ \&\ y=1)\vee (\top\in x^*\ \& \ y=0).\] By the assumption,
 $(\top\in x^{**}) \vee (\top \in x^*)$ for all $x\in
\Pow(\{\top\})$, so that trivially \[(\forall x \in
\Pow(\{\top\}))(\exists y\in \{0,1\}) (x,y)\in F.\] By GUP-CZF, this
gives \[(\exists y\in \{0,1\})(\forall x \in \Pow(\{\top\}))
(x,y)\in F,\] which yields a contradiction (consider $x=\{\top\},
x=\emptyset$,  for $y=0$, $y=1$, respectively).

Of course, this implies that GUP is inconsistent with the principle
of excluded middle EM in CTT (see also \cite{Petit}), in CZF with excluded middle
for restricted formulas, or equivalently, with
\[(\forall p \in \Pow(\{\top\}))p\cup p^*=\{\top\}\ \ \  (\textbf{[R]EM})\]

\noindent (again, considered in IZF, [R]EM is equivalent
to the full law of excluded middle).

In the following, we shall write CZF$^*$ or CTT$^*$ for a generic fixed extension of CZF or CTT, respectively,  that is compatible with the generalised uniformity
principle.
For simplicity, we call  any  extension of this kind a `constructive
setting' (this terminology is quite improper, given that CZF$^*$ may  be taken to be given by CZF plus the impredicative unbounded separation scheme Sep).
With `intuitionistic setting' we indicate any of CTT$^*$, CZF$^*$, IZF, HHA.

We shall make free use of the fact that all the settings
that we consider may consistently be extended with the negation of the
(restricted) De Morgan law.  The
assertion that a space of a certain type cannot be proved to have a
certain property in a certain setting will invariably  be proved by showing that in the setting
extended with some compatible non-classical principle (as GUP, or
$\neg$[R]DML), the assumption that the space has the
property is contradictory.

For definiteness, in what follows we always argue in the setting of constructive set theory.
The proof for a different system for which a given result is claimed, is obtained by the expected modifications of the given argument.

\section{{\normalsize Boolean and De Morgan locales/formal spaces}}\label{BFS}
By exploiting the generalised uniformity principle it is shown in \cite{CSSV} that there is a formal space that CTT cannot prove to be set-presented;
in \cite{G1} it is shown by other means that  the system CZF cannot prove the so-called  `double-negation'  formal space $\Pow(\{\top\})_{\neg\neg}$ (see Section \ref{OtherEx}) to be set-presented\footnote
{\label{grayson}
As  noted in \cite{G1},
R. Grayson \cite{Gra} had obtained a corresponding result for certain formulations of intuitionistic set theory without the powerset axiom.}.
We shall see in Section \ref{OtherEx} that $\Pow(\{\top\})_{\neg\neg}$ and the space considered in \cite{CSSV} are in fact isomorphic, and that their associated frame is a complete Boolean algebra. The same argument given  in \cite{G1} can then be used to show that any formal space whose associated frame is a non-trivial (complete) Boolean algebra
cannot be proved to be set-presented over the basic set of axioms of CZF\footnote{This observation is essentially due to S. Vickers; in fact also the corresponding of this  result was  known to Grayson \cite{Gra} in connection with the  set theories he considered (cf. footnote \ref{grayson}).}.

In this section, using the generalised uniformity principle  we give a simple proof  that no system  CZF$^*$, CTT$^*$
can prove a non-trivial Boolean formal space to be set-presented. In particular, thus, this holds for CZF$^*$=CZF+REA+PA+Sep. A similar result is  also shown to hold  for De Morgan (or extremally disconnected) formal spaces.

Further independence results concerning
overtness, compactness and existence of points are also  obtained. Aside from Theorem \ref{DeMorgannon-compact2}, these make no use of
the consistency of GUP with the given setting and  hold true,
\emph{mutatis mutandis}, also with respect to topos logic (HHA), or IZF. All
results in which the generalised uniformity principle is involved, which thus
only concern the constructive settings, will be marked with GUP.

Call a formal space $S$ such that $Sat(S)$ is a Boolean frame a
\textit{Boolean formal space}. From now on, let for simplicity $S$
have a top basic element $1_S$.  If $S$ is set-presented, also the
enlargement of its base with a top  element $1_S$ can be proved to be set-presented
(this is  proved in type theory using type-theoretic choice \cite{CuSC},
in constructive set theory exploiting the Subset Collection scheme).
Thus, a Boolean formal space is one such that
$1_S=_S U\cup U^*$ for all $U\in \Pow(S)$.

A formal space $S$ is \textit{De Morgan} if it satisfies
$1_S\cov U^{**}\cup U^*$ for all $U\in \Pow(S)$. Classically, a topological space is extremally disconnected iff its  frame of open subsets is De Morgan \cite{J82}.
Obviously, $S$ Boolean implies $S$ De Morgan.

For $p\in \Pow(\{\top\})$, we shall suggestively write $P$ to stand
for $\top\in p$, while `$\forall p$' will always stand for `for all $p$ in
$\Pow(\{\top\})$'. We set
\[U_P=\{x\in S: x=1_S\ \& \ P\}\equiv \{x\in S: x=1_S\ \&\ \top\in p\}.\]

\noindent
Note that,
in any formal space, $\{1_S\}^*= S^*=_S \emptyset$, and $\emptyset
^*= S=_S \{1_S\}$.

\begin{lem} Let $S$ be any formal space.

\begin{itemize}
\item [$i.$] \label{lemmaBoole} \textit{If $(\forall p)(\exists x)x\in U_P\cup U_P^*\ \& \ \neg (x\cov \emptyset)$, then  $(\forall p)P\vee \neg P$, i.e., $(\forall
p\in \Pow(\{\top\})) p\cup p^*=\{\top\}$.}

  \item [$ii.$] \label{lemmaDeMorgan}
\textit{If $(\forall p)(\exists x)x\in U^{**}_P\cup U_P^*\ \& \ \neg (x\cov \emptyset)$, then  $(\forall p)\neg\neg P\vee \neg P$, i.e., $(\forall p\in \Pow(\{\top\}))$ $p^{**}\cup p^*=\{\top\}$.}

\end{itemize}

\end{lem}

\proof
From $x\in U_P$ one gets $P$. From $x\in U_P^*$ and $\neg (x\cov \emptyset)$ one obtains $\neg P$ as follows: assuming $P$, one has  $U_P=\{x\in S: x=1_S\}$, so that $U_P^*\cov \emptyset$; together with  $x\in U_P^*$ and $\neg (x\cov \emptyset)$, this yields a contradiction, so that $\neg P$. Finally, assuming $\neg P$ gives $U_P=\emptyset$, and thus also $U^{**}_P=_S \emptyset $; by $x\in U^{**}_P$, $\neg (x\cov \emptyset)$ one derives $\neg\neg P$. The reader may then easily fill in the details.
\enproof

In \cite{FS79} one  finds an `arrow-theoretic' proof that no
Boolean frame may have points, unless classical logic is accepted.
Here is another formulation of that proof, and the corresponding fact for De Morgan locales.

\begin{prop}\label{point} No
 Boolean formal space can have a point unless [R]EM is accepted in
CZF$^*$ (EM in CTT$^*$, HHA, IZF). No De Morgan formal space can have a point unless [R]DML is accepted in CZF$^*$ (DML in CTT$^*$, HHA, IZF).
\end{prop}

\proof Assume Boolean formal space $S$ has a point.

\begin{itemize}

\item[] $S$ Boolean implies  $(\forall p)1_S \cov U_P \cup U^*_P$;

\item[] $S$ has a point $\alpha$ implies that  $(\forall p)(\exists a)a \in U_P \cup U^*_P\ \& \ a\in \alpha$.
\end{itemize}

\noindent By $a\in \alpha$ one has $\neg (a\cov \emptyset)$, and then
one concludes using Lemma \ref{lemmaBoole}. The proof for the De Morgan case is similar.
\enproof

As  the settings in consideration (CZF$^*$, CTT$^*$,  HHA, IZF)
can be extended consistently by  $\neg$[R]DML ($\neg$DML),  in these settings De Morgan
locales/formal spaces cannot be proved to have points. This implies that
no such formal space is a topological space, i.e., no non-trivial De Morgan (in particular Boolean)  frame can be obtained as the frame of opens of a non-empty (inhabited) topological space. Classically, of course,
the lattice of open subsets of any discrete (non-empty) space is a Boolean
frame with points.

Despite this result, at least in HHA/IZF, there are Boolean locales
that are \textit{proper}, i.e., such that, for all $U\in \Pow(S)$,
$S\cov U$ implies $\exists a\in U$ (see \cite{FS79}). Properness is
a stronger formulation of non-triviality.

A formal space $(S,\cov)$ is \textit{open} (or overt, or has a
positivity predicate \cite{J84, Sa87, CSSV}) iff there is a predicate
$\Pos(x)$, for $x$ in $S$, satisfying
\begin{itemize}
  \item [$i.$] $\Pos (a)$ and  $a \cov U$ imply $(\exists b \in S)b\in  U \
\& \ \Pos (b) $ \hspace{2.3cm} (\textit{monotonicity});
  \item [$ii.$] $a \cov U$ implies $a \cov U^+ \equiv\{b\in U: \Pos (b)\}
\hspace{3.8cm}  $ (\textit{positivity}).
\end{itemize}

\noindent(Note that classically, all formal spaces are open, with
$\Pos(a)\equiv \neg (a\cov \emptyset)$).
Then, although a Boolean locale can be proper, it cannot be open.

\begin{prop}\label{open}
No non-trivial De Morgan formal space $S$ (in particular, no non-trivial Boolean formal space) can be proved to be open in the intuitionistic settings  considered.
\end{prop}

\proof Assume $S$ is open. Then, by positivity,
$1_S\cov\{1_S\}^+$. Assume  $1_S\in \{1_S\}^+$, so that $\Pos(1_S)$
holds. Then,

\begin{itemize}

\item[] $S$  De Morgan implies  $(\forall p)1_S \cov U^{**}_P \cup U^*_P$;

\item[] by $\Pos(1_S)$ and monotonicity of $\Pos$,  for all $p$ there is  $a \in U^{**}_P \cup U^*_P$ with $\Pos(a)$.
\end{itemize}

Thus, since $\Pos(a)$ implies $\neg (a\cov \emptyset)$ (by
monotonicity), by Lemma \ref{lemmaBoole} one obtains $(\forall p)\neg\neg P\vee \neg P$.
As the settings under consideration can be extended
with the negation of (restricted) De Morgan law, one
has that $1_S\in \{1_S\}^+$ leads to a contradiction in the extended
setting, so that $\{1_S\}^+=\emptyset$. But this cannot be, as, by
positivity $1_S\cov \{1_S\}^+$, and we assumed the space to be
non-trivial.
\enproof

This proof shows that Boolean formal spaces can be proper only because the
elements one extracts from each cover of the whole space are  not
required to be  different from the empty open, let alone positive.

\begin{thm}[GUP]\label{set-presented}\label{DeMorgannon-sp} Let $S$ be
\begin{itemize}
  \item [$i.$] a non-trivial Boolean formal space, or
  \item [$ii.$] a non-trivial De Morgan formal space  such that $(\forall x\in
S) x\cov \emptyset \vee \neg (x\cov \emptyset);$
\end{itemize}
then $S$  cannot be proved to be set-presented (inductively generated) in CZF$^*$  or in CTT$^*$.
\end{thm}

\proof $i.$ Assume $S$ has a set-presentation in CZF$^*$.  Then

\begin{itemize}

\item[] $S$ Boolean implies:  $(\forall p)1_S \cov U_P \cup U^*_P$;

\item[] $S$ set-presented implies: $(\forall p)(\exists i\in
I(1_S))C(1_S,i) \subseteq U_P \cup U^*_P$;

\item[] assuming GUP, this implies:  $(\exists i\in I(1_S))(\forall p)C(1_S,i) \subseteq U_P \cup U^*_P$.

\end{itemize}

\noindent In particular, taking $p=\{\top\}$ so that $P\equiv\top \in \{\top\}$ is true, this gives
\[(*)\ \ \ \forall x \in C(1_S,i)(x=1_S \vee x\cov \emptyset).\]

\noindent Then assume $x=1_S\in C(1_S,i)$. Since $C(1_S,i) \subseteq
U_P \cup U^*_P$, and since $S$ is non-trivial, by Lemma
\ref{lemmaBoole}, one gets (in CZF$^*$+GUP) $(\forall p)P\vee \neg P$. We saw that
GUP is incompatible with [R]EM. Therefore, by $(*)$, one gets $\forall x\in C(1_S,i)x\cov
\emptyset$; but $1_S \cov C(1_S,i)$, and we assumed $S$ to be
non-trivial, whence $S$ is not set-presented in CZF$^*$+GUP. This shows that $S$ cannot be proved to be set-presented in CZF$^*$.

$ii.$ Assume $S$ is set-presented in CZF$^*$. Note first that if $S$ satisfies $(\forall x\in
S) x\cov \emptyset \vee \neg (x\cov \emptyset)$, but has no top element $1_S$, also
the isomorphic formal space with a top \cite{CuSC} will satisfy the given decidability condition, as, for the top element, we have by hypothesis that $\neg (1_S\cov \emptyset)$.
Then, in CZF$^*$+GUP, one
finds $i\in I(1_S)$ such that $(\forall p)C(1_S,i) \subseteq
U^{**}_P \cup U^*_P$. Let $x\in C(1_S,i)$, and assume $\neg (x\cov
\emptyset)$. By Lemma \ref{lemmaDeMorgan}, one obtains $(\forall p)\neg\neg
P\vee \neg P$. As this contradicts GUP, one has $\neg
\neg (x\cov \emptyset)$. By the decidability of $x\cov \emptyset$, it
follows that $x\cov \emptyset$ for all $x\in C(1_S,i)$. But
$S$ is non-trivial, so that it (is not set-presented in CZF$^*$+GUP and thus) cannot be proved to be set-presented in CZF$^*$.
\enproof

\noindent
An example of a  De Morgan non-Boolean formal space satisfying the condition in Theorem \ref{set-presented} is presented in the next section.

This theorem shows that the formal spaces that cannot be inductively
generated consist not just of few pathological cases. In particular, one has:

\begin{cor}\label{dnn}
Given any non-trivial formal space $S$, the formal space
$S^{**}\equiv(S, $ $\cov^{**})$, where $a\cov ^{**}U \iff \{a\}^{**}\cov
U^{**}$, can not be proved set-presented in CTT$^*$,
CZF$^*$.
\end{cor}

\noindent The space $S^{**}$ is indeed the Boolean formal subspace
corresponding to the ${**}-$nucleus on the frame defined by
$S$, i.e., the space associated with the frame of `regular' elements of $Sat(S)$ (see e.g. \cite{FS79,J82} for a discussion of the ${**}-$nucleus on
a locale $L$).\footnote{In the literature on locales this nucleus is also
known as the Booleanization of $S$. There are many Boolean
sublocales of a given locale $L$, but each of them can be seen as
defined by a ${**}-$nucleus over a closed sublocale of $L$
\cite{FS79}.} $S^{**}$ is non-trivial as soon as $S$ is non-trivial.

These subspaces/nuclei are, also classically, a peculiarity of
locale theory (as opposed to point-set topology), since, given
\textit{any} locale $L$, the ${**}-$nucleus on $L$ yields the least
dense sublocale of $L$. (This need not exist in a topological
space; consider e.g. the real line: the rationals and the
irrationals define dense disjoint subspaces).

\begin{remark}
Although not set-presentable, Boolean formal spaces are constructively
useful: an example of the use of a Boolean formal space to obtain a
concrete (constructive and predicative) description of ideal
non-effective objects can be found in
\cite[Theorem 6.1]{CL06}.
\end{remark}

The following proposition shows that being Boolean also conflicts with being compact.

\begin{prop}\label{compact}
No non-trivial Boolean formal space $S$  can be compact unless [R]EM in CZF$^*$, EM in CTT$^*$, HHA,  IZF, is accepted.
\end{prop}

\proof Assume $S$ is compact.

\begin{itemize}

\item[] $S$ Boolean implies  $(\forall p)1_S \cov U_P \cup U^*_P$;

\item[] $S$ compact implies that, for all $p$, there is a finite $u_0$ such that $1_S \cov u_0\subseteq U_P \cup
U^*_P$.
\end{itemize}

\noindent  It is a standard fact that $u_0\subseteq V \cup W$, with
$u_0$ finite, implies intuitionistically $u_0=v_0\cup w_0$, with
$v_0\subseteq V$, $w_0\subseteq W$ both finite \cite{C92}. Thus we
have finite $v_0\subseteq U_P, w_0\subseteq U^*_P$, with $1_S\cov
v_0\cup w_0$; moreover, `finite' implies  `either empty or
inhabited'. By cases: $v_0, w_0=\emptyset$ cannot be, by
non-triviality. Then one of the following alternatives holds:

\begin{itemize}

\item[$1.$] $v_0, w_0$ inhabited, or

\item[$2.$] $v_0$ inhabited and $w_0=\emptyset$, or

\item[$3.$] $v_0=\emptyset$ and $ w_0$ inhabited.

\end{itemize}

\noindent The first and second case directly give $P\vee \neg P$.
For the last, assuming $P$ one gets $w_0\cov \emptyset$, that
together with $v_0=\emptyset$, gives $1_S\cov \emptyset$, so that,
by non-triviality of $S$, $\neg P$, and then again $P\vee \neg P$. Therefore, if $S$ is compact,
the law $(\forall p)P\vee \neg P$ holds.
\enproof

In CZF$^*$, or CTT$^*$, more generally, no non-trivial Boolean space $S$ can be proved to be locally
compact,  since locally compact formal spaces are set-presented
\cite{A}, and by Theorem \ref{set-presented} no non-trivial Boolean space $S$ can be proved to be
set-presented in these settings (`more generally': any Boolean $S$ is regular, and a
compact regular locale is locally compact, e.g. \cite{J82}).
Classically (e.g. in ZF), every finite discrete space has a compact
Boolean frame of opens.

\medskip

So far the generalised uniformity principle has only been used to show that a constructive system cannot prove that formal spaces of a certain type can be set-presented. We conclude this section with two other important consequences of  the consistency of this principle with  the constructive settings we are considering.

First let us note that, contrary to what one may expect, \textit{in HHA or IZF, De Morgan
locales can be compact}: the classical result (due to M. Stone) that
the (compact) frame $Idl(B)$ of ideals over a Boolean algebra $B$ is
De Morgan if and only if $B$ is complete (e.g. \cite{Johnstone80}) is topos-valid.
The following is one half of this result, formulated for
formal spaces. Recall that in e.g.  IZF, frames
and set-generated class-frames come to the same thing, so that
$Sat(S)$ is carried by a set for every space $S$.

\begin{prop}[Stone]\label{idealicBaDM}
In any of the intuitionistic settings we are considering, let  $S$ be a Boolean formal space. In the context of CZF$^*$, or CTT$^*$, assume further that $Sat(S)$ is (carried by) a set. Then the formal space $S_\beta\equiv
(Sat(S),\cov_\beta)$, with $U\cov_\beta \{U_i\}_{i\in I} \iff U\cov
U_{i_1}\cup ... \cup U_{i_n}$ for $\{i_1,...,i_n\}$ a (possibly
empty) finite subset of $I$, is a compact De Morgan formal space.
\end{prop}

\proof The proof that $S_\beta$ is a formal space
is left to the reader. One has to prove \[S\cov_\beta (\{U_i\}_{i\in
I})^{*_\beta}\cup (\{U_i\}_{i\in I})^{*_\beta*_\beta}\] for any given set $\{U_i\}_{i\in
I}$ of elements of $Sat(S)$. Routine calculations show that
\begin{equation}
(\{U_i\}_{i\in I})^{*_\beta}=_\beta \{(\bigcup_{i\in I} U_i)^{*_S}\}
\end{equation}
($^{*_S}$ is pseudo-complementation in $S$), so that
\begin{equation}
((\{U_i\}_{i\in I})^{*_\beta})^{*_\beta}=_\beta \{(\bigcup_{i\in I} U_i)^{*_S*_S}\}.
\end{equation}
Since $S$ is Boolean, $(\bigcup_{i\in I} U_i)^{*_S*_S}=_S
\bigcup_{i\in I} U_i$. As \[S \cov (\bigcup_{i\in I} U_i)^{*_S}\cup
(\bigcup_{i\in I} U_i)=_S (\bigcup_{i\in I} U_i)^{*_S}\cup
(\bigcup_{i\in I} U_i)^{*_S*_S},\] one has $S\cov_\beta
\{(\bigcup_{i\in I} U_i)^{*_S}, (\bigcup_{i\in I} U_i)^{*_S*_S}\}$
by definition of $\cov_\beta$ (note that pseudo-complements are
saturated). Therefore, by 1,2 above, $S\cov_\beta (\{U_i\}_{i\in
I})^{*_\beta} $ $\cup (\{U_i\}_{i\in I})^{*_\beta *_\beta}$, as wished.
\enproof

\noindent Despite this fact, one has:

\begin{prop}\label{DeMorgannon-compact1}
No non-trivial De Morgan formal space $S$ such that $(\forall x\in
S)$ $ \neg\neg (x\cov \emptyset) \vee \neg (x\cov \emptyset)$ is compact unless [R]DML holds in CZF$^*$ (DML holds in CTT$^*$, HHA, IZF).
\end{prop}

\proof  For all $p\in \Pow(\{\top\})$, one
finds a finite $v_0$ with $v_0\subseteq U^{**}_P\cup U^*_P$,
$1_S\cov v_0$. One has, in particular, $(\forall x\in v_0) \neg\neg
(x\cov \emptyset) \vee \neg (x\cov \emptyset)$. By a general
intuitionistic principle (see \cite[Lemma 2.4]{J84}), this gives
$(\forall x\in v_0) (\neg\neg (x\cov \emptyset)) \vee (\exists x\in
v_0)(\neg (x\cov \emptyset))$. Since $v_0$ is finite, we get
$\neg\neg(\forall x\in v_0) (x\cov \emptyset) \vee (\exists x\in
v_0)(\neg (x\cov \emptyset))$. It cannot be that $x\cov \emptyset$ for
all $x\in v_0$. Thus,  there is $x\in v_0$ with $\neg (x\cov
\emptyset)$, so that one concludes by Lemma \ref{lemmaDeMorgan}.
\enproof

As a consequence, in a topos that does not satisfy De Morgan law,
no frame of the form $Idl(B)$, with
$B$ complete Boolean algebra, can have a base satisfying the
decidability condition in the above proposition. Note also that such
frames  are examples of De Morgan frames that are never intuitionistically Boolean, given that no
Boolean frame can be proved compact.

Using the generalised uniformity principle, the above proposition may be
strengthened.

\begin{thm}[GUP]\label{DeMorgannon-compact2}
No non-trivial De Morgan formal space $S$  can be proved to be
compact in  CZF$^*$, CTT$^*$.
\end{thm}

\proof Using GUP, one has that a finite subset
$u_0=\{x_1,...,x_n\}$  of $S$ exists such that
$(\forall p)1_S\cov u_0\subseteq U_P^*\cup U_P^{**}$ ($u_0$ is non-empty, as
$S$ is non-trivial). Assume $\neg (x_1\cov \emptyset) \vee...\vee \neg
(x_n\cov \emptyset)$. By Lemma \ref{lemmaDeMorgan}, one has that
$(\forall p)\neg\neg P\vee \neg P$ holds. We saw that this
principle is incompatible with GUP, so that $\neg(\neg (x_1\cov
\emptyset) \vee...\vee \neg (x_n\cov \emptyset))$. This gives $\neg
\neg (x_1\cov \emptyset)\ \&...\&\ \neg\neg (x_n\cov \emptyset)$, that
is $\neg \neg(x_1\cov \emptyset\ \&...\&\ x_n\cov \emptyset)$. On
the other hand, from $1_S\cov u_0$ one gets $\neg(x_1\cov \emptyset\
\&...\& \ x_n\cov \emptyset)$, so that $S$ (is not  compact in  CZF$^*$+GUP, and hence) cannot be proved to be compact in CZF$^*$.
\enproof

\noindent Recall that the \emph{Gleason cover} of a compact regular formal space $S$ is a minimal surjection $\gamma S\to S$, with $\gamma S$ a compact, regular, De Morgan formal space  \cite{Johnstone80,J82}. It then follows from  Theorem \ref{DeMorgannon-compact2} that, in contrast with what happens in a topos, \emph{for no non-trivial  compact regular formal space $S$ the Gleason cover of $S$ can be constructed in CZF$^*$, CTT$^*$.}

By  Theorem \ref{DeMorgannon-compact2}  it also follows  that no non-trivial  frame can be
assumed to be carried by a set in a constructive setting (see \cite{CuExistenceSC} for a more direct proof).

\begin{cor}[GUP]\label{Noframeset2}
Every non-trivial frame $Sat(S)$, for $S$ Boolean, is carried by a
proper class in  CZF$^*$, CTT$^*$. Thus, no non-trivial
frame $Sat(S)$ may be proved to have a set of elements in these
contexts.
\end{cor}

\proof If the collection of elements of $Sat(S)$ could be proved to be
constructively a set, by Proposition \ref{idealicBaDM}   the formal
space $S_\beta$ would be compact and De Morgan, contradicting
Theorem \ref{DeMorgannon-compact2}. Now assume a frame $Sat(S)$ has
a set of elements; then all frames $Sat(S')$, for $S'$ a subspace of
$S$, are carried by a set, too, so that also $Sat(S^{**})$ should
be.
\enproof

\begin{remark}
 The property of Boolean formal spaces that has been exploited in the proofs of Propositions \ref{point},  \ref{compact}, and Theorem \ref{set-presented}, is that the whole space $S$ is covered by $U_P \cup U^*_P$, for all $p$ in $\Pow(\{\top\})$ (this is also true for Proposition \ref{open}, if one proves the result just for the Boolean case). It is easy to check that a morphism $f:S\longrightarrow S'$, with $S$ Boolean, preserves pseudocomplements. It follows that whenever such a morphism exists, one also has  $(\forall p)1_{S'} \cov V_P \cup V^*_P$, with $V_P\equiv \{x\in S': x=1_{S'}\ \& \ P\}$. Then, Propositions \ref{point}, \ref{open}, \ref{compact}, and Theorem \ref{set-presented} hold true more generally if one replaces the Boolean space $S$  with any non-trivial codomain of a morphism with Boolean domain.
Similar considerations also hold in connection with the results concerning De Morgan spaces in Propositions \ref{point}, \ref{open}, \ref{DeMorgannon-compact1}, and Theorems \ref{set-presented}, \ref{DeMorgannon-compact2}, when $f:S\longrightarrow S'$ is any morphism that preserves  pseudocomplements  (in particular, when $f$ defines an \emph{open} continuous functions  of locales/formal spaces \cite{J84}).
\end{remark}

In contrast with the Boolean case, one cannot hope to prove that every subspace of a De Morgan formal space is De Morgan: classically, an extremally disconnected space may have Hausdorff
subspaces that are not extremally disconnected.
In \cite{J79}, the following law is considered: for
all propositions $P,Q$
\[(P\to Q)\vee (Q\to P).\]
This principle is stronger than De Morgan's (take $Q$ to be $\neg
P$), and is inherited by the internal logic of sheaf subtoposes
\cite{J79}. Call \emph{strongly De Morgan} a formal space such that
the associated frame models this formula. A strongly De Morgan
formal space is De Morgan. It is easy to prove that the class of
strongly De Morgan formal spaces is closed for subspaces.

\section{{\normalsize Dedekind--MacNeille completions}} \label{OtherEx}

Given a set $S$, and any (class-)relation $R(a, U)$, for $a\in S$ and $U\in \Pow(S)$,
one may define $R$ to be set-presented precisely as for coverings. Let  $\varPhi(P)$ be an instance of a law  in one  variable $P$ that is incompatible with GUP, e.g. $\varPhi(P)\equiv P \vee \neg P$ (in CZF$^*$,  $p\cup p^*=\{\top\}$, for $p\in \Pow(\{\top\})$).

\begin{prop}[GUP]\label{relations}
Let $S$ be a set, and let $R\subseteq S\times \Pow(S)$ be such that, for some $a$ in $S$, $\neg R(a,\emptyset)$, and  $(\forall p)R(a,U_{\varPhi(P)})$, where $U_{\varPhi(P)}\equiv \{x\in S: x=a \ \& \ \varPhi(P)\}$. Then $R$ cannot be proved to be set-presented in CZF$^*$, CTT$^*$.
\end{prop}

\proof Assume $R$ is set-presented by $C(x,i)$ with $x\in S$ and $i\in I(x)$. By GUP, there is $i\in I(a)$ such that $(\forall p)C(a,i)\subseteq U_{\varPhi(P)}$ and $R(a, C(a,i))$. Assume $x\in C(a,i)$. Then $(\forall p)\varPhi(P)$. By hypothesis this  contradicts GUP. Therefore,  $C(a,i)=\emptyset$, and $R(a, \emptyset)$, against what we have assumed.
\enproof

Recall that the Dedekind--MacNeille completion of a partial order makes it possible to embed a given partially ordered set in a complete lattice preserving meets and joins that exist (see e.g. \cite{Sa89,TvD}).
Given a partially ordered set $(S,\leq)$, one may define a relation $R_c(x,U)$  by letting \[R_{c}(x, U)\iff (\forall
y)[(\forall u\in U) u\leq y]\to x\leq y\iff x\in \bigcap_{U\subseteq
\downarrow y} \downarrow y.\]
To have that $R_{c}$ is a covering relation, the \emph{Dedekind--MacNeille covering}, the partial order has to satisfy some further conditions. In particular, if $S$ is a  Heyting algebra this is always the case. The  frame $Sat(S,R_c)$ of saturated subsets of the formal topology $(S,R_c)$ is then the complete lattice in which the Heyting algebra $S$ is embedded via $e:S\longrightarrow Sat(S,R_c)$, $e(a)={\cal S}(\{a\})$. Recall that, as over any frame, an implication operation making $Sat(S,R_c)$ a complete Heyting algebra can be defined by letting $U\to V\equiv \{a\in S: a\downarrow U\cov V\}$. The  Heyting algebra structure of $Sat(S,R_c)$ then extends that of $S$ (see e.g. \cite{TvD}, vol. II).

T. Coquand has suggested\footnote{On the occasion of the presentation of the material in the preceding sections at the workshop ``Trends in constructive mathematics", Chiemsee
(Germany) June 19-23, 2006.}
that no Dedekind--MacNeille covering can be constructively proved to be
set-presented (see also \cite{CSSS}).
We prove here that  this holds for every relation $R_c(x,U)$ on a given poset $(S,\leq)$, but with a further hypothesis.

\begin{prop}[GUP]\label{GDM}
Let $(S,\leq)$ be a partial order having at least one element $a$ that is not the least of $S$, and that is `stable', in the sense that $\neg \neg (a \leq x)$ implies $a\leq x$, for all $x$. Then the relation $R_{c}(x,U)$ cannot be proved to be set-presented in CZF$^*$, CTT$^*$.
\end{prop}

\proof  By Proposition \ref{relations}, it suffices to show that $\neg R_{c}(a,\emptyset)$ and $(\forall p)R_{c}(a,U_{\varPhi(P)})$, with  $U_{\varPhi(P)}\equiv \{a:\varPhi(P)\}$, and $\varPhi(P)\equiv P \vee \neg P$. If $R_{c}(a,\emptyset)$, then $a\in \bigcap_{\emptyset\subseteq
\downarrow y} \downarrow y.$ This gives  $(\forall y\in S) a\leq y$, against the hypothesis.
For the second, let $U_{\varPhi(P)}\subseteq \downarrow y$, and assume $\neg (a \leq y)$ and $\varPhi(P)$. Then $U_{\varPhi(P)}= \{a\}\subseteq \downarrow y$, so that $a\leq y$, against what we have assumed. This gives $\neg \varPhi(P)$.
As $\neg\neg \varPhi(P)$ is intuitionistically provable, we get $\neg \neg (a \leq y)$, whence $a\leq y$.  We conclude that $a\in \bigcap_{U_{\varPhi(P)}\subseteq
\downarrow y} \downarrow y$, for all $p$, i.e., $(\forall p)R_{c}(a, U_{\varPhi(P)})$.
\enproof

Note that this proof is, in essence, a simplification and a generalization of the proof for the special case considered in \cite{CSSV}.

\begin{cor}[GUP]
If $(S,\leq)$ is any  poset with at least two elements and a decidable order relation, then $R_{c}(x,U)$ cannot be proved to be set-presented in  CZF$^*$, CTT$^*$.
\end{cor}

These results may be used to produce examples of non-De Morgan formal spaces that cannot be constructively set-presented.

\begin{cor}[GUP]
Let $H\equiv (S,\wedge, \vee, \to, 0, 1)$ be a Heyting algebra satisfying the hypothesis in Proposition \ref{GDM} (w.r.t. the partial order associated with $H$). Assume in $H$  De Morgan law is false, i.e.,  there is $b\in S$ such that $1\neq b^*\vee b^{**}$. Then the relation $R_{c}(x,U)$, defining the Dedekind--MacNeille cover on $H$, defines a non-De Morgan formal space that cannot be proved to be set-presented in  CZF$^*$, CTT$^*$.
\end{cor}

\proof  As already recalled, the (set-generated) frame $Sat(S,R_c)$ associated with the Dedekind--MacNeille cover defined over an Heyting algebra $H$ is a complete Heyting algebra in which the Heyting algebra operations extend the corresponding operations on $H$. As $b^*=b\to 0$ one can conclude.
\enproof

The set $T=\{0,\frac{1}{2}, 1\}$ endowed with the natural order is a non-Boolean  Heyting algebra. The Dedekind--MacNeille cover over this poset defines a De Morgan non-Boolean formal space $(T,\cov_{DM})$. That $T$ is non-Boolean again follows by the fact that the complete Heyting algebra $Sat(T)$ is such that the Heyting algebra operations are extensions of the corresponding operations of $T$. To prove that $(T,\cov_{DM})$ is De Morgan, i.e., that  for every $U\in \Pow(T)$,  $1\in \bigcap K(U)$, with $K(U)\equiv \{\downarrow y: U^*\cup U^{**}\subseteq \downarrow y\}$, first one notes that $U^*\cup U^{**}\subseteq \downarrow 1$;
unwinding the definitions, one then proves
that assuming $1\not \in U^*\cup U^{**}$ leads to a contradiction, so that $\neg 1\not \in U^*\cup U^{**}$ (sketch:  $1\not \in U^*\cup U^{**}$ implies $\neg (1\in U^*) \ \& \  \neg (1\in U^{**})$. The  second conjunct yields $\neg (\frac{1}{2}\not \in U^*\ \&\ 1\not \in U^*)$, that in turn gives $(\frac{1}{2}\not \in U\ \&\ 1\not \in U)$; on the other, by the first conjunct one gets  $\neg (\frac{1}{2}\not \in U\ \&\ 1\not \in U)$, so that a contradiction is reached). Therefore, $\downarrow \frac{1}{2}\not \in K(U)$ and $\downarrow 0\not \in K(U)$, whence $1\in \bigcap K(U)$. By Proposition \ref{GDM} (or Theorem \ref{DeMorgannon-sp}) one has that $(T,\cov_{DM})$ is not constructively set-presentable.

\begin{remark}
If in the hypotheses of Proposition \ref{relations},  $R$ is a covering $\cov$ on a set $S$,  and if  $a\equiv 1_S=_S S$, then for every morphism $f:S\longrightarrow S'$, by $1_s\cov U_{\varPhi(P)}$ one gets  $1_{S'}\cov U'_{\varPhi(P)}$, with $U'_{\varPhi(P)}\equiv \{x\in S': x=1_{S'}\ \&\ \varPhi(P)\}$. Therefore, if $S'$ is non-trivial, by Proposition \ref{relations} it is not   set-presentable. As an immediate corollary one has in particular that no formal space defined by the Dedekind--MacNeille cover on a poset $(S,\leq)$ with the properties in Proposition \ref{GDM}, and such that $a$ is also the greatest element of $(S,\leq)$, may have points (as points are in a bijective correspondence with morphisms from the given topology to the set-presentable topology $\Pow(\{\top\})$), and that every non-trivial formal subspace of $S$ is not set-presentable (if $S'\equiv (S,\cov')$ is a subspace of $S$, letting $e(x)=\{x\}$ for all $x\in S$  defines a morphism $e:S\longrightarrow S'$).
\end{remark}

In fact, the next result can be obtained  without any reference to the uniformity principle (see also \cite{CSSS}).

\begin{prop}\label{}
Let $S$ be a Dedekind--MacNeille topology defined on a poset with the property in Proposition \ref{GDM}, and having $a$ as greatest element. If $S$ has a point, or is  compact, then [R]EM in CZF$^*$ (EM in CTT$^*$, HHA, IZF) holds. Furthermore,  $S$ cannot be proved to be open in the  intuitionistic settings considered.
\end{prop}

\proof Follows immediately by the fact that it holds $\neg R_{c}(a,\emptyset)$ and  $R_{c}(a,U_{\varPhi(P)})$ for all $p$, with  $\varPhi(P)\equiv(p\cup p^*=\{\top\})$ (cf. proof of  Proposition \ref{GDM}).
\enproof

We conclude this section with a discussion of Open Problem 4.5
of \cite{CSSV}. Taking  $S=\{0,1\}$ with the natural order, and $a\cov_{DM} U\equiv R_{c}(a, U)$ we obtain  the formal space that is in \cite{CSSV} shown not to be  set-presentable in CTT. Let us denote this space by $S_{DM}$.

We already pointed out that a  uniform  method for the definition of  products of arbitrary formal spaces is generally regarded as being beyond constructive means. Open Problem 4.5
of \cite{CSSV}  asked  whether at least  the particular product of $S_{DM}$ with itself is predicatively definable. The answer in this (indeed very special) case is yes. We show this with a slight detour.
The `double negation' formal topology
$\Pow(\{\top\})_{\neg\neg}$, is  defined by
\[S=\{\top\}, \top\cov U\iff \neg\neg (\top\in U).\]
In \cite{G,G1} it is shown that the system CZF cannot prove $\Pow(\{\top\})_{\neg\neg}$ to be set-presentable.
 Note that $\Pow(\{\top\})_{\neg\neg}$ is isomorphic with the Boolean formal space  $\Pow
(\{\top\})^{**}$ of regular elements of $\Pow(\{\top\})$ (cf. section \ref{BFS}).

It is then easy to see that $\Pow(\{\top\})_{\neg\neg}$ and
$S_{DM}$ are the `same' formal space.

\begin{lem}\label{DM=NN}
$S_{DM}\cong \Pow(\{\top\})_{\neg\neg}$.
\end{lem}
\proof It is an  exercise in  intuitionistic
logic to prove that $0=_{S_{DM}}\emptyset$, $\neg \neg (1\in U) \iff
1\cov_{DM} U$, and $\neg \neg (\top \in U)$ implies $1\cov _{DM} \{1:
\top \in U\}$. It follows that the homomorphisms $f:S_{DM}\longrightarrow \Pow(\{\top\})_{\neg\neg}$ and $g:\Pow(\{\top\})_{\neg\neg}\longrightarrow S_{DM}$, given by $f(0)=\emptyset$, $f(1)=\{\top\}$, and $g(\top)=
\{1\}$, yield the required isomorphism.
\enproof

Now  the product
$\Pow(\{\top\})_{\neg\neg}\ \times\ \Pow(\{\top \})_{\neg\neg}$ in
FSp is simply $\Pow(\{\top\})_{\neg\neg}$ itself:
$\Pow(\{\top\})_{\neg\neg}$ is `almost' a terminal object in FSp, if a morphism with $\Pow(\{\top\})_{\neg\neg}$ as domain exists, then it is unique.

\section*{{\normalsize Conclusion}}

The generalised uniformity principle has been systematically exploited  in this paper to obtain
non-derivability results for the main formal systems for constructive mathematics,  in particular with the aim of  distinguishing  topos-valid from intuitionistic generalised predicative mathematics.
It has already been pointed out that it is  improper to define `constructive' an extension of CZF or CTT that is consistent with the generalised uniformity principle, as e.g. CZF plus the highly impredicative unbounded separation principle is one such extension.
On the other hand, it may be reasonable to define
`non-constructive' a result that cannot be derived within some extension of CZF or CTT that is compatible with GUP.
This note has thus shown that some standard topos-valid results are in fact non-constructive in this sense.
A further important  result, valid in any topos,  that turns out to be a non-constructive theorem in the present sense is described in \cite{CuExistenceSC}.

\section*{{\normalsize Appendix: The generalised uniformity principle in type theory}}

In the type-theoretic context the generalised uniformity principle reads informally as follows: given any set $I$, and any mapping $R$ into the type of propositions $PROP$ taking a proposition and an element of $I$ as arguments,
if a mapping $F$ is given from $PROP$ to  $I$ such that $R(P,F(P))$ holds for all $P$, then one may find an element $\bar i\in I$ such that $R(P,\bar i)$ holds for all $P$. We denote this version of the principle by \textbf{GUP-CTT}.
Observe that, due to the propositions-as-sets identification, the type $PROP$ of propositions
may be replaced in GUP-CTT by the type $SET$ of sets. In \cite{CSSV}
this principle is formulated implicitly and it is claimed that
GUP-CTT can be `added' consistently to type theory.
Models, due to T. Coquand, of  type theory  validating this form of the uniformity principle have then been discussed in \cite{Petit}.

Formally, GUP-CTT can
be expressed in type theory (more specifically, in the logical
framework \cite{NPS90, CDPS}) by the addition of two constants UP1,
UP2  as follows (cf. \cite{Petit}):

\medskip
\noindent $UP1:(I: SET)\to R:(PROP\to El(I)\to PROP)\to F:(PROP\to El(I))\to$\\
$G:(P:PROP\to El(R(P,F(P))))\to El(I);$

\medskip

\noindent
$UP2:(I: SET)\to R:(PROP\to El(I)\to PROP)\to F:(PROP\to El(I))\to$\\
$G:(P:PROP)\to El(R(P,F(P)))\to ((P:PROP)\to El(R(P,UP1(IRFG)))).$

\medskip

\noindent \textbf{Acknowledgements.}
I began the research described in
this note while  visiting  the university of Birmingham in 2005.
Thanks go to the ``Fondazione Gini" for  financially  supporting
that visit.
In that occasion I enjoyed stimulating conversations with Steven
Vickers on this and related topics.

I am grateful to Thierry Coquand, Milly Maietti, and Pino Rosolini
for a helpful and interesting e-mail discussion on the consistency
of the generalised uniformity principle with type theory, and to Thomas
Streicher for letting  me know about the proof of the consistency of
this principle with various extensions of CZF in Benno van den Berg's thesis.
Thanks also go to Peter Aczel and Paul Taylor for their useful comments on a previous
draft of this note.


{\small
\begin{thebibliography}{10}


\bibitem{A}
P. Aczel, Aspects of general topology in constructive set
theory, Ann. Pure Appl. Logic \textbf{137}  1--3,
3--29 (2006).

\bibitem{AczelCuri}
P. Aczel and G. Curi, On the $T_1$  axiom and other separation
properties in constructive point-free and point-set topology, Ann. Pure Appl. Logic, 161 (2010), pp. 560-569.

\bibitem{AR} P. Aczel and M. Rathjen,  Notes on Constructive Set Theory,
 Mittag-Leffler Technical Report No.40, 2000/2001.

\bibitem{BennoThesis}
B. van den Berg, Predicative topos theory and models for
constructive set theory, PhD thesis, University of Utrecht, 2006.

\bibitem{BergMoerdijk}
B. van den Berg and I. Moerdijk, Aspects of predicative algebraic set
theory II: realizability, Theor. Comp. Science, to appear. [Available from: http://arxiv.org/abs/0801.2305].

\bibitem{C92}
T. Coquand, An intuitionistic proof of Tychonoff theorem,
J. Symb. Log.  \textbf{57}  1,  28--32 (1992).

\bibitem{CL06}
T. Coquand and H. Lombardi, A logical approach to abstract algebra. A
survey, Math. Struct. in Comp. Science \textbf{16}, 885--900 (2006).

\bibitem{CDPS}
T. Coquand, P. Dybjer, E. Palmgren and A. Setzer, Type-theoretic
foundation of constructive mathematics. In preparation.


\bibitem{CSSS}
T. Coquand,  S. Sadocco, G. Sambin, J. Smith,  Formal
topologies on the set of first-order formulae, J. Symb. Log. \textbf{65} 3, 1183--1192 (2000).


\bibitem{CSSV}
T. Coquand,  G. Sambin, J. Smith and  S. Valentini, Inductively
generated formal topologies, Ann. Pure Appl. Logic \textbf{124}  1--3,  71--106 (2003).



\bibitem{CuSC}
G. Curi, Exact approximations to Stone-\v{C}ech
compactification, Ann. Pure Appl. Logic \textbf{146} 2--3,
 103--123 (2007).

\bibitem{CuExistenceSC}
G. Curi, On the existence of Stone-\v{C}ech compactification, J. Symb. Log., to appear.


\bibitem{FG82} M. Fourman and R. Grayson, Formal Spaces, in:  The
L.E.J. Brouwer Centenary Symposium, edited by: A.S. Troelstra and D. van
Dalen (North Holland, 1982), pp. 107--122.


\bibitem{FS79} M. Fourman and D.S. Scott, Sheaves and logic, in:
Applications of sheaves, edited by: M. Fourman et al. (Springer
LNM \textbf{753}, Springer-Verlag, 1979), pp. 302--401.

\bibitem{G} N. Gambino, Sheaf interpretations for generalised
predicative intuitionistic systems, Ph.D. thesis, University of
Manchester, 2002.

\bibitem{G1} N. Gambino, Heyting-valued interpretations for
Constructive Set Theory,  Ann. Pure Appl. Logic \textbf{137} 1--3,
 164--188 (2006).

\bibitem{Gra} R.J. Grayson, Forcing in intuitionistic systems
without powerset, J. Symb. Log.  \textbf{48}, 670--682 (1983).


\bibitem{J79}
P. T. Johnstone, Conditions related to De Morgan's law, in:
Applications of sheaves, edited by: M. Fourman et al. (Springer
LNM \textbf{753}, Springer-Verlag, 1979),  pp. 479--491.

\bibitem{Johnstone80}
P. T. Johnstone, The Gleason cover of a topos, I. J. Pure Appl. Algebra 19, 171-192 (1980).


\bibitem{J82}
P. T. Johnstone,  Stone Spaces (Cambridge University Press,
1982).

\bibitem{J83}
P. T. Johnstone, The points of pointless topology,
Bull. Am. Math. Soc. \textbf{8} 1, 41--53 (1983).


\bibitem{J84}
P. T. Johnstone,  Open locales and exponentiation,
Contemp. Math. \textbf{30},  84--116  (1984).


\bibitem{ML84}
P. Martin-L\"of, Intuitionistic Type Theory.
 Studies in Proof Theory (Bibliopolis, Napoli,
1984).


\bibitem{Pa} E. Palmgren, Predicativity problems in point-free
topology. In:   Proceedings of
the Annual European Summer Meeting of the Association for Symbolic
Logic, held in Helsinki, Finland, August 14--20, 2003,  Lecture
Notes in Logic \textbf{24}, ASL, edited by: V. Stoltenberg-Hansen et al. (AK Peters Ltd, 2006), pp. 221-231.

\bibitem{Lubarsky2006}
R.S. Lubarsky, CZF and Second Order Arithmetic, Ann. Pure Appl. Logic \textbf{141}  1--2, 29--34 (2006).


\bibitem{MM}
S. MacLane and I. Moerdijk, Sheaves in Geometry and Logic - A
First Introduction to Topos Theory (Springer, 1992).

\bibitem{NPS90}
B. N\"ordstrom, K. Peterson, J. Smith, {\it Programming in
Martin-L\"{o}f's Type Theory} (Clarendon Press), Oxford, 1990.

\bibitem{Petit}
B. Petit, Polymorphic type theory and Uniformity Principle.
Preprint.

\bibitem{Rathjen06}
M. Rathjen, Realizability for constructive Zermelo-Fraenkel set theory. In:   Proceedings of
the Annual European Summer Meeting of the Association for Symbolic
Logic, held in Helsinki, Finland, August 14--20, 2003,  Lecture
Notes in Logic \textbf{24}, ASL, edited by: V. Stoltenberg-Hansen et al. (AK Peters Ltd, 2006), pp. 282-314.

\bibitem{Ro}
G. Rosolini, About modest sets, Int. J. Found. Comp.
Sci \textbf{1}, 341--353 (1990).

\bibitem{Sa87}
G. Sambin, Intuitionistic formal spaces - a first
communication, in:  Mathematical Logic and its Applications,
edited by: D. Skordev (Plenum,  1987), pp. 187--204.


\bibitem{Sa89}
G. Sambin, Pretopologies and completeness proofs, J. Symb. Log. \textbf{60} 3, 861--878 (1995).

\bibitem{Streicher}
T. Streicher, Realizability models for CZF+ $\neg$ Pow, unpublished note.

\bibitem{TvD}
A. Troelstra and D. van Dalen, Constructivism in mathematics,
an introduction, Volumes I,II. Studies in logic and the foundation of
mathematics \textbf{121, 123} (Amsterdam etc., North-Holland, 1988).


\end{thebibliography}
}
\end{document}